\documentclass[a4paper]{amsart}
\usepackage[all]{xypic}
\newtheorem{thm}{Theorem}[section]
\newtheorem{cor}[thm]{Corollary}
\newtheorem{lem}[thm]{Lemma}
\newtheorem{prop}[thm]{Proposition}
\theoremstyle{definition}

\theoremstyle{remark}

\numberwithin{equation}{section}

\newcommand{\Z}{{\mathbb{Z}}}
\newcommand{\Q}{{\mathbb{Q}}}

\newcommand{\Tor}{\mathop{\textrm{\rm Tor}}}
\newcommand{\Hom}{\mathop{\textrm{\rm Hom}}}
\newcommand{\Ext}{\mathop{\textrm{\rm Ext}}}
\begin{document}

\title[Homological Localizations of Eilenberg-Mac\,Lane Spectra]
{Homological Localizations \\ of Eilenberg-Mac\;Lane Spectra}%
\author{Javier J.~Guti\'errez}\thanks{The author was supported by
MEC-FEDER grant MTM2004-03629}
\address{Departament d'\`Algebra i Geometria \\ Universitat
de Barcelona \\ Gran Via, 585 \\ \newline E-08007 Barcelona,
Spain}
\email{javier.gutierrez@ub.edu}

\begin{abstract}
We discuss the Bousfield localization $L_E X$ for any spectrum $E$
and any $HR$-module $X$, where $R$ is a ring with unit. Due to the
splitting property of $HR$-modules, it is enough to study the localization
of Eilenberg--Mac\,Lane spectra. Using general results about stable
$f$-localizations, we give a method to compute the localization of an
Eilenberg--Mac\,Lane spectrum $L_E HG$ for any spectrum $E$ and any abelian
group $G$. We describe $L_E HG$ explicitly when $G$ is one of the following:
finitely generated abelian groups, $p$-adic integers, Pr\"ufer groups,
and subrings of the rationals. The results depend basically on
the $E$-acyclicity patterns of the spectrum $H\Q$ and the spectrum
$H\Z/p$ for each prime~$p$.
\end{abstract}

\maketitle
\section{Introduction}

Homological localizations were first defined by Adams \cite{Ada73}. Bousfield developed
the theory further by proving the existence of homological localizations in the category of
spaces \cite{Bou75} and in the category of spectra \cite{Bou79b}.

Given any spectrum $E$, a homological localization functor with respect to
$E$ is a homotopy idempotent transformation $L_E\colon Ho^s\longrightarrow Ho^s$,
where $Ho^s$ is the stable homotopy category, that turns $E$-homology equivalences
into homotopy equivalences in a universal way. Homological localizations
are special cases of $f$-loca\-li\-za\-tions in the sense of \cite{Dro96} and commute
with the suspension operator.

In \cite{CG05}, we presented a general study of $f$-localizations of $HR$-module spectra
and discussed the preservation of several structures under the effect of these functors. In
this paper, we restrict our attention to homological localizations in order to obtain more
explicit results. In fact, we translate to spectra some of the results of \cite{Bou82}, by using
ideas of \cite{CG05} to simplify the arguments.

In \cite{Bou79b}, Bousfield determined the homological localizations of
connective spectra with respect to connective homology theories.
A spectrum $X$ is connective if $\pi_k (X)=0$ for $k<0$.
If either $E$ or $X$ fail to be connective, then $L_E X$ is
somehow unpredictable. For example, the spectrum $L_K S$, where $K$
denotes complex $K$-theory and $S$ is the sphere
spectrum, has infinitely many nonzero homotopy groups in both
positive an negative dimensions (see \cite[Theorem 8.10]{Rav84} or
\cite[Corollary 5.15]{CG05}).

We study $L_E X$ where $E$ is any homology theory (not
necessarily connective) and
$X$ is any $HR$-module spectrum for a ring $R$ with unit. Since any $HR$-module
splits as a wedge of suspensions of Eilenberg--Mac\,Lane spectra, we focus
on the study of $L_E HG$ for any homology theory $E$ and any abelian group
$G$, where $HG$ denotes the Eilenberg--Mac\,Lane spectrum associated
to $G$. We describe all possible homological localizations in the
case of finitely generated abelian groups and other groups, including
the $p$-adic integers, the Pr\"ufer groups $\Z/p^{\infty}$, and subrings of the
rationals. For example, in the case of the spectrum $H\Z$, by the general approach
of \cite{CG05} we know any of its localization has at most one nonzero
homotopy group and that this group has the structure of a \emph{rigid ring} in
the sense of \cite{CRT00}. We prove
that, for homological localizations of $H\Z$, the only rigid rings that appear are subrings
of the rationals or products of $p$-adic integers for diferent primes.

The computations of these localizations depend on the $E$-acyclicity patterns of the
spectra $H\Z/p$ and $H\Q$, and on the set of primes $p$ such that $G$ is uniquely
$p$-divisible, similarly as in \cite{Bou82}.

\bigskip
\noindent\textbf{Acknowledgements.} I am especially indebted to A.~K.  Bousfield for
ecouraging me to write this paper and for sharing his insight. I would also like
to thank Carles Casacuberta and Mark Hovey for many useful conversations.

\section{Homological localization of spectra}
We will work in the stable homotopy category of spectra $Ho^s$ (see \cite{Ada74}).
Any spectrum $E$ in $Ho^s$ gives rise to a homology theory defined as
$E_k(X)=\pi_k(E\wedge X)$ for any spectrum $X$ and any $k\in\Z$. Homological localization
with respect to the homology theory $E$ is a functor that transforms homology equivalences
with respect to this theory into homotopy equivalences in a universal way.
It is unique up to homotopy and idempotent.

A  map of spectra $f\colon X\longrightarrow Y$ is an
\emph{$E$-equivalence} if the map $f_*\colon E_k(X)\longrightarrow
E_k(Y)$ is an isomorphism for all $k\in\Z$. A spectrum $X\in Ho^s$
is called \emph{$E$-acyclic} if $E_k(X)=0$ for all $k\in\Z$, that
is, if $E\wedge X$ is contractible. A spectrum $Z$ is
\emph{$E$-local} if each $E$-equivalence $X\longrightarrow Y$
induces induces a homotopy equivalence $F(Y,Z)\simeq F(X,Z)$, or
equivalently if $F(W,Z)=0$ for each $E$-acyclic spectrum $W$,
where $F(X,Y)$ denotes the function spectrum from $X$ to $Y$. An
\emph{$E$-localization} of a spectrum X is a map $\eta_X\colon
X\longrightarrow L_E X$, where $X$ is an $E$-local spectrum and
$\eta_X$ is an $E$-equivalence. Homological localization is
universal in the following sense: the localization map $\eta_X$ is
initial among maps from $X$ to $E$-local spectra and it is
terminal among all $E$-equivalences with domain~$X$.

The class of all the $E$-acyclic spectra for a given spectrum $E$ is denoted by
$\langle E \rangle$ and called the \emph{Bousfield class} or the
\emph{acyclicity class} of $E$.
Given two spectra $E$ and $F$, the $E$-localization functor and the
$F$-localization functor are equivalent if and only if $\langle
E\rangle=\langle F\rangle$. By Ohkawa's theorem, there is only a set of Bousfield
classes \cite{DP01}, and therefore a set of non-equivalent homological localization
functors.

\section{Acyclicity patterns of $H\Z/p$ and localization of $HR$-modules}
In this section we study how the $E$-acyclicity
patterns of $H\Z/p$ determine the localization $L_{E\Z/p}X$ for any spectrum $E$
and any $HR$-module spectrum $X$.
For any spectrum $E$ and any abelian group $G$, let $EG=E\wedge MG$, where
$MG$ is the Moore spectrum associated to $G$.

A spectrum $E$ is called a \emph{stable $R$-GEM} if it is homotopy
equivalent to a wedge of suspensions of Eilenberg--Mac\,Lane spectra, i.e.,
$E\simeq \vee_{k\in\Z}\Sigma^k HA_k$, where each $A_k$ is an $R$-module
(hence, each $HA_k$ is an $HR$-module spectrum). If $R=\Z$, then stable
$\Z$-GEMs are called simply \emph{stable GEMs}. The Eilenberg--Mac\,Lane spectrum $HG$ is
an $R$-GEM if $G$ is an $R$-module. The stable $R$-GEMs are precisely the spectra that admit a module
structure over the ring spectrum $HR$ (see for example \cite[Proposition
4.4]{CG05}).

The splitting property of $HR$-modules allows us to
describe their localization easily. Note that every $HR$-module is an
$H\Z$-module trivially via the morphism $\Z\rightarrow R$ that sends
the unit of $\Z$ to the unit of the ring $R$. And also that homological
localizations commute with suspension, i.e., $L_E \Sigma^k
HG\simeq\Sigma^k L_E HG$ for all $k\in\Z$, since the desuspension of an $E$-equivalence
is again an $E$-equivalence \cite[Proposition 2.4]{CG05}.

\begin{prop}
$L_E(\vee_{k\in\Z}\Sigma^k HA_k)\simeq \vee_{k\in\Z}(\Sigma^k L_E
HA_k)$ for any spectrum $E$.
\label{comsus}
\end{prop}
\begin{proof}
The spectrum $\vee_{k\in\Z}(\Sigma^k L_E HA_k)$ is $E$-local, since
the natural map
$$
\bigvee_{k\in\Z}\Sigma^k L_E HA_k\longrightarrow \prod_{k\in\Z}\Sigma^k L_E HA_k
$$ is a homotopy equivalence, because by \cite[Theorem 5.6]{CG05} for each value of
$k$, at most two nonzero homotopy groups appear in $L_E HA_k$. Now, the map
$$
\bigvee_{k\in\Z}\Sigma^k HA_k\longrightarrow \bigvee_{k\in\Z} \Sigma^k L_E HA_k
$$
is an $E$-equivalence, because it is a wedge of $E$-equivalences.
\end{proof}

In \cite{Bou79a}, Bousfield showed that
\begin{equation}
\langle E \rangle=\langle
E\Q\rangle \vee \bigvee_{p\in\mathcal{P}}\langle E\Z/p\rangle,
\label{boudec}
\end{equation}
for any spectrum $E$, where $\mathcal{P}$ is the set of all primes.
In fact, what this means essentially is that we can recover $L_E X$
for any $E$ and $X$ from information on what happens rationally,
$L_{E\Q}X$, and at each prime, $L_{E\Z/p}X$.
The following result of Bousfield \cite[Proposition 2.9]{Bou79b} ilustrates this fact.
Recall that a commutative diagram
of spectra
$$
\xymatrix{
X \ar[r]^f  \ar[d]_h & Y \ar[d]^g \\
Z \ar[r]_i & W
}
$$
is an \emph{arithmetic square} if there is a map $j\colon
W\longrightarrow \Sigma X$ such that
$$
\xymatrix{
X\ar[r]^-{(f,h)} & Y\wedge Z \ar[r]^-{(g,-i)} & W \ar[r]^j & \Sigma
X
}
$$
is a cofiber sequence of spectra.
\begin{prop}
For all spectra $E$ and $X$, there is an arithmetic square
\begin{equation}
\xymatrix{L_E X \ar[r] \ar[d]&
\prod_{p\in\mathcal{P}} L_{E\Z/p} X \ar[d] \\ L_{E\Q}X \ar[r] &
L_{E\Q}(\prod_{p\in\mathcal{P}} L_{E\Z/p} X),
}
\label{arithsq}
\end{equation}
where $\mathcal{P}$ is the set of all primes. $\hfill\qed$
\end{prop}

The $E\Q$-localizations were completely determined in
\cite{Bou79b}. For any spectrum $E$, all these localizations are
equivalent to rationalization. In fact, $L_{E\Q}X=L_{M\Q}X=X\wedge M\Q$ for
all $E$ and $X$.

The computation of $L_{E\Z/p}X$ for any $E$ and any $HR$-module spectrum
$X$ depends on the $E$-acyclicity types of the spectrum $H\Z/p$ for each prime $p$.
Note that if $\Sigma^i H\Z/p$ is $E$-acyclic for some $i\in\Z$, then $\Sigma^k H\Z/p$ is $E$-acyclic
for all $k\in\Z$, since homological localizations commute with suspension.

\begin{prop}
If $H\Z/p$ is $E$-acyclic, then $L_{E\Z/p}X=0$ for any $HR$-module spectrum $X$.
\label{prop01}
\end{prop}
\begin{proof}
It is enough to check that $M\Z/p$ is $E$-acyclic, because in this case
$E\Z/p\wedge X\simeq E\wedge M\Z/p\wedge X=0$. The spectrum $M\Z/p\wedge X$ is obviously
$E\Q$-acyclic and $E\Z/q$-acyclic for $q\ne p$. In the case $q=p$, we have that
$$
M\Z/p \wedge X \wedge E\Z/p\simeq X'\wedge
H\Z\wedge M\Z/p\wedge E\Z/p=0
$$
since $X$ is an $H\Z$-module and therefore splits as $H\Z\wedge X'$ for some spectrum $X'$, and
$H\Z\wedge E\Z/p\simeq E\wedge H\Z/p=0$. Now using the decomposition (\ref{boudec}), we have
that $X$ is $E\Z/p$-acyclic.
\end{proof}

\begin{lem}
If $H\Z/p$ is not $E$-acyclic and $f\colon X\longrightarrow Y$ is an
$E\Z/p$-equivalence, then it is an $H\Z/p$-equivalence.
\label{lem02}
\end{lem}
\begin{proof}
Since homological localizations commute with suspension, if we smash $f$ with any spectrum
the resulting map is an $E\Z/p$-equivalence. In particular, if we smash with the spectrum
$H\Z$, the map $f\wedge H\Z$ induces an equivalence
$$
E\Z/p \wedge H\Z\wedge X  \simeq E\Z/p \wedge H\Z \wedge Y.
$$
The spectrum $E\Z/p\wedge H\Z\simeq E\wedge H\Z/p$ is an $H\Z/p$-module,
so its homotopy gropus are $\Z/p$-vector spaces and it splits as a wedge $\vee_{k\in
I}\Sigma^k H\Z/p$ (there may be repetitions in the index set $I$) and this wedge
is non-trivial since by hypothesis $E\wedge H\Z/p\ne 0$. Hence, $f$ induces an equivalence
$$
\bigvee_{k\in I}\Sigma^k H\Z/p\wedge X\simeq \bigvee_{k\in I}\Sigma^k H\Z/p\wedge
Y
$$
which turns $f$ into an $H\Z/p$-equivalence.
\end{proof}
The following theorem allows us to compute the localization $L_{E\Z/p}$ of connective spectra or $HR$-modules when
the spectrum $H\Z/p$ is not $E$-acyclic.

\begin{thm}
If $H\Z/p$ is not $E$-acyclic, then $L_{E\Z/p}X\simeq L_{M\Z/p}X$ for every
spectrum $X$ that is connective or an $HR$-module.
\label{thm01}
\end{thm}
\begin{proof}
If $X$ is a connective spectrum, then $L_{M\Z/p}X\simeq L_{H\Z/p}X$ (see \cite[Theorem 3.1]{Bou79b}).
The localization map $X\longrightarrow L_{M\Z/p} X\simeq
L_{H\Z/p}X$ is an $M\Z/p$-equivalence and therefore and $E\Z/p$-equivalence. Moreover,
the spectrum $L_{H\Z/p}X$ is $E\Z/p$-local since by Lemma
\ref{lem02} every $H\Z/p$-local spectrum is $E\Z/p$-local.

If $X$ is an $HR$-module, the result follows from the above and
Proposition \ref{comsus}.
\end{proof}

The case $L_{M\Z/p}X$ can be computed using \cite[Proposition 2.5]{Bou79b}:
\begin{prop}
For any spectrum $X$, we have that
$$
L_{M\Z/p}X\simeq F(\Sigma^{-1}M\Z/p^{\infty}, X),
$$
and there is a splittable exact sequence
$$
0\longrightarrow \Ext(\Z/p^{\infty},\pi_k (X))\longrightarrow
\pi_k(L_{M\Z/p} X)\longrightarrow
\Hom(\Z/p^{\infty},\pi_{k-1}(X))\longrightarrow 0
$$
for any $k\in\Z$. $\hfill\qed$
\end{prop}
In the particular case when $X$ is an Eilenberg--Mac\,Lane spectrum $HG$, we have that
$L_{M\Z/p}HG\simeq HA\vee \Sigma HB$ where $A\cong\Ext(\Z/p^{\infty},G)$
and $B\cong\Hom(\Z/p^{\infty},G)$.

\section{Localizations of Eilenberg-Mac\,Lane spectra}

In the study of homological localizations of $HR$-module spectra,
we can focus our attention on the particular case of homological
localizations of Eilenberg--Mac\,Lane spectra $L_E HG$, by
Proposition \ref{comsus}.  Homological localizations are a particular
example of homotopical localizations or $f$-localizations. These
localizations in the stable homotopy category have been studied in \cite{CG05}.
In that paper, we proved that the localization of any Eilenberg--Mac\,Lane spectrum has
at most two nonzero homotopy groups in dimensions zero and one (see \cite[Theorem 5.6]{CG05}).
Thanks to the Bousfield arithmetic square, to compute $L_E HG$ it is enough to
determine $L_{E\Z/p}HG$ for every prime $p$, since in the rational case
$L_{E\Q}HG=H(\Q\otimes G)$ for any spectrum $E$.

An abelian group $G$ is called \emph{uniquely $p$-divisible} if for every
$g\in G$ there exists a unique $h\in G$ such that $g=ph$. This condition is
equivalent to saying that $\Z/p \otimes G=0$ and $\Tor(\Z/p,G)=0$.

\begin{lem}
For any spectrum $E$, the group $\pi_k(E)$ is uniquely $p$-divisible for all $k\in\Z$
if and only if $E\Z/p=0$.
\label{lem01}
\end{lem}
\begin{proof}
The result follows using the exact sequence
$$
\Z/p\otimes \pi_k(E)\longrightarrow \pi_k(E\Z/p)\longrightarrow
\Tor(\Z/p, \pi_{k-1}(E)),
$$
which is valid for every $k\in\Z$.
\end{proof}

As a particular case, we have that the abelian group $(H\Z)_k (E)$
is uniquely $p$-divisible for all $k\in\Z$ if and only if $H\Z/p$
is $E$-acyclic. Note also that if $\pi_k(E)$ is uniquely
$p$-divisible, then $(H\Z)_k(E)$ is uniquely $p$-divisible.

\begin{prop}
If $H\Z/p$ is not $E$-acyclic, then $L_{E\Z/p}HG=0$ if and only if $G$
is uniquely $p$-divisible.
\label{prop02}
\end{prop}
\begin{proof}
If $L_{E\Z/p}HG=0$, then $E\wedge H\Z/p\wedge MG=0$. Since
$E\wedge H\Z/p$ is an $H\Z/p$-module spectrum, we have that
$$
E\wedge H\Z/p \wedge MG=\vee_{k\in I}\Sigma^k H\Z/p\wedge MG=0.
$$
Therefore, $H\Z/p\wedge MG=M\Z/p\wedge HG=0$ and thus $G$ is uniquely $p$-divisible by
Lemma \ref{lem01}.

On the other hand, if $G$ is uniquely $p$-divisible, then by Lemma \ref{lem01}
we have that $HG\wedge M\Z/p=0$ and hence $L_{E\Z/p}HG=0$.
\end{proof}

By means of Theorem \ref{thm01} and Proposition \ref{prop02}, one can now compute the localization
$L_{E\Z/p}HG$ depending on the $E$-acylicity patterns of $H\Z/p$. If $H\Z/p$ is $E$-acyclic, then
$L_{E\Z/p}HG=0$. If $H\Z/p$ is not $E$-acyclic, then $L_{E\Z/p}HG\simeq L_{M\Z/p}HG$ if $G$ is not
uniquely $p$-divisible and zero otherwise. The arithmetic square (\ref{arithsq}) in the case $X=HG$ is
the following:

\begin{equation}
\xymatrix{ L_E HG \ar[r] \ar[d]&
\prod_{p\in\mathcal{P}} L_{E\Z/p} HG \ar[d] \\ H(G\otimes \Q) \ar[r] &
M\Q \wedge(\prod_{p\in\mathcal{P}} L_{E\Z/p} HG),
}
\label{arithsq2}
\end{equation}
where $\mathcal{P}$ is the set of all primes $p$ such that $H\Z/p$
is not $E$-acyclic and $G$ is not uniquely $p$-divisible.

\begin{thm}
Let $A_p=\Ext(\Z/p^{\infty},G)$, $B_p=\Hom(\Z/p^{\infty},G)$, and let $\mathcal{P}$ be the set
of primes such that $H\Z/p$ is not $E$-acyclic and $G$ is not uniquely
$p$-divisible. For any spectrum $E$ and any abelian group $G$, we have the following:
\begin{itemize}
\item[(i)] If $H\Q$ is $E$-acyclic, then
$$
L_E HG=\prod_{p\in\mathcal{P}}(HA_p\vee \Sigma
HB_p).
$$
\item[(ii)] If $H\Q$ is not $E$-acyclic, then there is a cofiber sequence of spectra

$$
L_E HG \longrightarrow H(\Q\otimes G)\vee \prod_{p\in\mathcal{P}}(HA_p\vee \Sigma HB_p)
\longrightarrow M\Q\wedge \prod_{p\in\mathcal{P}}(HA_p\vee \Sigma HB_p).$$
\end{itemize}
\label{thm02}
\end{thm}
\begin{proof}
The result follows from Proposition \ref{prop01}, Proposition \ref{prop02}, Theorem \ref{thm01} and the arithmetic square
(\ref{arithsq2}).
\end{proof}

\section{Some examples}
In this section, we compute homological localizations of Eilenberg--Mac\,Lane spectra and
$HR$-module spectra in some concrete examples.
First, we compute $L_E X$ for some non-connective homology theories $E$
and any $HR$-module spectrum~$X$.
\subsection{Localization with respect to $n$-th Morava $K$-theory
$K(n)$} Let $n\ge 0$, $p$ a fixed prime and $K(n)$ the spectrum of the $n$-th Morava $K$-theory at $p$.
Recall that $\pi_* K(n)\cong \Z/p\,[v^{-1}_n, v_n]$ where
$|v_n|=2(p^n-1)$ for $n\ge 1$.

If $n=0$, then $K(0)=H\Q=M\Q$ and so $L_{K(0)}HG=H(G\otimes\Q)$.
In the case $n\ge 1$ we know that $K(n)\wedge M\Q=0$ and $H\Z/p$ is $K(n)$-acyclic for
every prime $p$, because $K(n)\wedge H\Z/p=0$ for all primes $p$ (see for example [Rav84, Theorem 2.1]).
Thus $L_{K(n)}HG=0$ for $n\ge 1$. Hence,

\begin{prop}
For any $HR$-module $X$, its localization with respect to $K(n)$ is
either zero if $n\ge 1$, or rationalization if $n=0$, i.e.,
$L_{K(0)}X=X\wedge M\Q$. $\hfill\qed$
\end{prop}

\subsection{Localization with respect to Johnson--Wilson spectra $E(n)$}
The Bousfield class of $E(n)$ splits as a wedge of Morava $K$-theories,
$\langle E(n)\rangle=\langle K(0)\vee\ldots\vee K(n)\rangle$ (see [Rav84, Theorem 2.1]);
therefore $L_{E(n)}HG=L_{K(0)}HG=H(G\otimes\Q)$ since $L_{K(i)}HG=0$ if $i\ge 0$.

\subsection{Localization with respect to complex $K$-theory}
The spectrum $H\Z/p$ is $K$-acyclic for every prime $p$ and $K\Q\ne 0$,
so $L_K HG=H(G\otimes \Q)$. Therefore, we infer the following:

\begin{prop}
For any $HR$-module spectrum $X$, its localization with respect to
$E(n)$ or $K$-theory is rationalization. $\hfill\qed$
\end{prop}

In the next examples, we use Theorem \ref{thm02} to compute all
the possible homological localizations of the spectrum $HG$ with
respect to any $E$ for some families of abelian groups. Given any
spectrum $E$ and any abelian group $G$, we have the following
acyclicity patterns that determine the localization $L_E HG$
completely. These patterns are the stable analogues of Condition I
and Condition II of \cite[Section 4]{Bou82}:

\begin{itemize}
\item Pattern I: $E\Q=0$ and $E\wedge H\Z/p=0$ for all primes $p$.
\item Pattern II: $E\Q\ne 0$ and $E\wedge H\Z/p=0$ for all primes $p$.
\item Pattern III: $E\Q=0$ and $E\wedge H\Z/p\ne 0$ for all primes $p$ in a set of primes~$\mathcal{P}$.
\item Pattern IV: $E\Q\ne 0$ and $E\wedge H\Z/p\ne 0$ for all primes $p$ in a set of primes~$\mathcal{P}$.
\end{itemize}
Note that if Pattern I holds, we have that $L_E HG=0$ for any abelian group $G$.
\subsection{Localizations of $H\Z$}
\label{infcic}
The abelian group of the integers is not uniquely $p$-divisible for any
prime $p$. If Pattern II holds, then $L_E H\Z=H\Q$. We have that
$$
\Hom(\Z/p^{\infty},\Z)=0\quad\mbox{ and }\quad
\Ext(\Z/p^{\infty},\Z)=\widehat{\Z}_p,
$$
where $\widehat{\Z}_p$ is the ring of $p$-adic integers.
If Pattern III holds, then $L_E
H\Z=H(\prod_{p\in\mathcal{P}}\widehat{\Z}_p)$. And if Pattern IV holds,
then taking $\pi_0$ in the square (\ref{arithsq2}) we have the following
pullback diagram of abelian groups:
$$
\xymatrix{ \pi_0(L_E H\Z)
\ar[r] \ar[d]& \prod_{p\in\mathcal{P}} \widehat{\Z}_p \ar[d] \\
\Q \ar[r] & \Q \otimes \prod_{p\in\mathcal{P}} \widehat{\Z}_p\,,}
$$
where $\mathcal{P}$ is the set of all primes $p$ such that $E\wedge H\Z/p\ne 0$.
So $L_E H\Z=H\Z_{\mathcal{P}}$.

In \cite[Theorem 5.12]{CG05} we proved that every $f$-localization of the
spectrum $H\Z$ has at most one nonzero homotopy group, which aquires the structure of
a rigid ring in the sense of \cite{CRT00}. A ring $A$ with unit is \emph{rigid} if
evaluation at $1$ induces an isomorphism of abelian groups $\Hom(A,A)\cong A$.
In the special case of homological localizations we get the following:

\begin{prop}
For any spectrum $E$, we have that $L_E H\Z$ is either zero or
$HA$, where the rigid ring $A$ is a subring of $\Q$ or a product
of $p$-adic integers for different primes. $\hfill\qed$
\end{prop}

\subsection{Localizations of $H\Z/p^k$ for a prime $p$}
\label{cic}
The group $\Z/p^k$ is uniquely $q$-divisible for every $q\ne p$ and
moreover $L_{E\Q}H\Z/p^k\simeq H(\Q\otimes \Z/p^k)=0$ for
all $p$. We have that
$$
\Hom(\Z/p^{\infty},\Z/p^k)=0\quad\mbox{ and }\quad \Ext(\Z/p^{\infty},\Z/p^k)=\Z/p^k,
$$
hence $L_E H\Z/p^k=0$ under Pattern II and $L_E H\Z/p^k=H\Z/p^k$ under
Pattern III or Pattern IV.

\subsection{Localizations of $H\Q$}
\label{rat}
The group $\Q$ is uniquely $p$-divisible for every prime $p$, so $L_{E\Z/p}H\Q=0$ for all $p$.
If Pattern III holds, then $L_E H\Q=0$ and $L_E H\Q=H\Q$ under Pattern II or Pattern IV.

\subsection{Localization of $H\Z_{\mathcal{R}}$ for a set of primes $\mathcal{R}$}
For every prime $p\in\mathcal{R}$, we have that
$$
\Hom(\Z/p^{\infty},
\Z_{\mathcal{R}})=0 \quad\mbox{ and }\quad \Ext(\Z/p^{\infty},
\Z_{\mathcal{R}})=\widehat{\Z}_p.
$$
In fact, $\Hom(\Z/p^{\infty}, G)=0$ if $G$ is a torsion-free
abelian group and $\Ext(\Z/p^{\infty},G)=0$ if and only if $G$ is
$p$-divisible. If Pattern II holds, then $L_E
H\Z_{\mathcal{R}}=H\Q$ because $\Q\otimes \Z_{\mathcal{R}}\cong
\Q$. If Pattern III holds, then $L_E
H\Z_{\mathcal{R}}=H(\prod_{p\in\mathcal{R}\cap\mathcal{P}}\widehat{\Z}_p)$.
And if Pattern IV holds, then $L_E
H\Z_{\mathcal{R}}=H\Z_{\mathcal{R}\cap\mathcal{P}}$, where
$\mathcal{P}$ is the set of all primes $p$ such that $H\Z/p$ is
not $E$-acyclic. Note that this case generalizes the cases of the
localization of $H\Z$ (when $\mathcal{R}=\emptyset$) and $H\Q$
(when $\mathcal{R}$ is the set of all primes).

\subsection{Localizations of $H\Z/p^{\infty}$}
\label{divg}
The group $\Z/p^{\infty}$ is uniquely $p$-divisible for every prime
$q\ne p$. In this case, $L_{E\Q}H\Z/p^{\infty}=0$ since $\Q\otimes \Z/p^{\infty}=0$ for all $p$.
We have that
$$
\Hom(\Z/p^{\infty},\Z/p^{\infty})=\widehat{\Z}_p
\quad\mbox{ and }\quad\Ext(\Z/p^{\infty},\Z/p^{\infty})=0.
$$
Thus, under Pattern II, $L_E H\Z/p^{\infty}=0$. If Pattern III holds, then
$L_E H\Z/p^{\infty}\simeq \Sigma H\widehat{\Z}_p$.
And if Pattern IV holds, then by Theorem \ref{thm02} we have a cofiber sequence of spectra
$$
L_E
H\Z/p^{\infty}\longrightarrow \Sigma H\widehat{\Z}_p\longrightarrow
\Sigma H\widehat{\Q}\longrightarrow \Sigma H(\widehat{\Q}/\widehat{\Z}_p),
$$
where $\widehat{\Q}\cong\widehat{\Z}_p\otimes\Q$ are the $p$-adic rationals.
Hence, $L_E H\Z/p^{\infty}=H(\widehat{\Q}/\widehat{\Z}_p)\simeq H\Z/p^{\infty}$.

\subsection{Localization of $H\widehat{\Z}_p$}
We only have to focus on the prime $p$, because $\widehat{\Z}_p$ is uniquely $q$ divisible
for all primes $q\ne p$. In this case, we have that
$$
\Hom(\Z/p^{\infty},\widehat{\Z}_p)=0\quad\mbox{ and
}\quad\Ext(\Z/p^{\infty},\widehat{\Z}_p)=\widehat{\Z}_p.
$$
If Pattern II holds, then $L_E H\widehat{\Z}_p=H\widehat{\Q}_p$.
And if Pattern III or Pattern IV hold, then $L_E H\widehat{\Z}_p=H\widehat{\Z}_p$.

\bigskip
The following table summarizes the results obtained for the homological localizations
of Eilenberg--Mac\,Lane spectra for different groups. The set $\mathcal{P}$ is the set of primes $p$
such that $H\Z/p$ is not $E$-acyclic.

$$
\renewcommand{\arraystretch}{1.5}
\begin{array}{|c|c|c|c|c|}
\hline
 & \mbox{Pattern I} & \mbox{Pattern II} & \mbox{Pattern III} &
 \mbox{Pattern IV} \\
\hline
L_E H\Z & 0 & H\Q & \prod_{p\in \mathcal{P}}H\widehat{\Z}_p &
H\Z_{\mathcal{P}} \\
\hline
L_E H\Z/p^k & 0 & 0 & H\Z/p^k & H\Z/p^k \\
\hline
L_E H\Q & 0 & H\Q & 0 & H\Q \\
\hline
L_E H\Z_{\mathcal{R}} & 0 & H\Q &
\prod_{p\in \mathcal{P}\cap \mathcal{R}}H\widehat{\Z}_p &
H\Z_{\mathcal{P}\cap \mathcal{R}} \\
\hline
L_E H\Z/p^{\infty} & 0 & 0 & \Sigma H\widehat{\Z}_p & H\Z/p^{\infty}
\\
\hline
L_E H\widehat{\Z}_p & 0 & H\widehat{\Q}_p & H\widehat{\Z}_p & H\widehat{\Z}_p \\
\hline
\end{array}
$$

\subsection{Localization of $HG$ where $G$ is a finitely
generated abelian group} Every finitely generated abelian group
splits as a direct sum $G=\oplus_{i=1}^n C_i$ where each $C_i$ is
either $\Z$ or $\Z/p^k$ for some prime $p$ and $k\ge 1$. Since
$HG\simeq\vee_{i=1}^n HC_i$, then $L_E HG=\vee_{i=1}^n L_E HC_i$
and the localization of each $HC_i$ is determined using the
results of sections \ref{infcic} and \ref{cic}.

\subsection{Localization of $HG$ where $G$ is a divisible abelian
group} If $G$ is a divisible abelian group, then $G\cong R\oplus T$, where $R=\oplus_i \Q$ and
$T=\oplus_p(\oplus_{j_p}\Z/p^{\infty})$. In this case $L_E HG\simeq L_E HR\vee L_E HT$.
Since $R$ is a retract of $\prod_i\Q$, we have that $L_E HR=HR$ or $L_E HR=0$ depending on
whether $H\Q$ is $E$-local or $E$-acyclic. The localization $L_E HT$ can be determined using
the exact sequence of abelian groups
$$
0\longrightarrow \Z_\mathcal{P}\longrightarrow \Q\longrightarrow \oplus_{p\in \mathcal{P}}\Z/p^{\infty}\longrightarrow 0
$$
together with the results of sections \ref{rat} and \ref{divg}, and the fact that homological
localizations preserve cofiber sequences.

\section{Localization of reduced Eilenberg-Mac\,Lane spectra}
In all the examples we have studied, except in the case of $H\Z/p^{\infty}$, all the
homological localizations of $HG$ have at most one nonzero homotopy group in dimension zero.
This property also holds when the group $G$ is abelian reduced. An abelian group is \emph{reduced}
if it does not have nontrivial divisible subgroups. We say that the Eilenberg--Mac\,Lane spectrum $HG$ is \emph{reduced}
if the group $G$ is reduced.

\begin{thm}
If $HG$ is reduced, then $L_E HG$ is either zero or $HA$ for some
abelian group $A$ and for any spectrum $E$.
\end{thm}
\begin{proof}
If $G$ is reduced, then $\Hom(\Z/p^{\infty}, G)=0$. The result follows now from Theorem \ref{thm02}.
\end{proof}

Any abelian group $G$ splits as a direct sum $G\cong G_1\oplus G_2$, where $G_1$ is the maximal
divisible subgroup of $G$ and $G_2$ is reduced. Morover, $G_1$ splits as a direct sum of $\Q$'s
and $\Z/p^{\infty}$ for several primes $p$. Hence, by Theorem \ref{thm02} and the results in the
previous section, the only possibility for the homological localization of an Eilenberg--Mac\,Lane
spectrum $HG$ to have a nonzero homotopy group in dimension one, is that some $\Z/p^{\infty}$ appears
as a factor of the decomposition of $G$ and that $L_E H\Z/p^{\infty}=\Sigma H\widehat{\Z}_p$.

\begin{cor}
If $\Z/p^{\infty}$ does not occur as a direct summand of $G$ for
any prime~$p$, then $L_E HG$ is either zero or $HA$ for some
abelian group $A$ and any spectrum $E$. $\hfill\qed$
\end{cor}


\end{document}